\documentclass{article}
\title{\bf $\mathcal{N}$-max infinite divisibility and $\mathcal{N}$-max stability} 
\author{Satheesh S\thanks{Corresponding author, Department of Applied Sciences,  Vidya Academy of Science and Technology, Thalakkottukara, Thrissur 680 501, India. email: \textit{ssatheesh1963@yahoo.co.in}, phone: +919447724700}\\
and\\
Sandhya E\thanks{Department of Statistics, Prajyoti Niketan College, Pudukad, Thrissur -680 301, India. email: \textit{esandhya@hotmail.com}}}

\def\lo{\left}
\def\ro{\right}
\def\pa{\paragraph}

\def\be{\begin{equation}}
\def\ee{\end{equation}}
\def\beq{\begin{eqnarray*}}
\def\eeq{\end{eqnarray*}}
\def\bea{\begin{eqnarray}}
\def\eea{\end{eqnarray}}

\begin{document}
\maketitle
\begin{abstract}  
Here we give a necessary and sufficient condition for the convergence to a random max infinitely divisible law from that of a random maximum. We then discuss random max-stable laws, their domain of max-attraction and the associated extremal processes.\\

\end{abstract}

\section{Introduction}
\setcounter{equation}{0}
\def\theequation{\thesection.\arabic{equation}} 

Balkema and Resnick (1977) had introduced the notion of max infinitely divisible (MID) laws and later Rachev and Resnick (1991) geometric max infinitely divisible (GMID) laws and geometric max stable laws (GMS), see also Mohan (1998). Since all distribution functions (\textit{d.f}) in $\mathbf{R}$ are MID, a discussion of MID laws is relevant for \textit{d.f}s in $\mathbf{R^d}, d\geq 2$, integer and the max operations are to be taken component wise. Thus in this paper, all \textit{d.f}s are assumed to be in $\mathbf{R^d}, d\geq 2$, integer. Rachev and Resnick (1991) also discussed certain connections between GMID/ GMS laws and extremal processes. From Balkema and Resnick (1977) we have the max-analogue of the classical de Finetti's theorem.

\pa{Theorem 1.1}A \textit{d.f} $F(x)$ is MID \textit{iff} for some \textit{d.f}s $\{G_n\}$ and constants $\{a_n>0\}$
\[F(x)=\lim_{n\to\infty}\exp\{-a_n(1-G_n(x))\}.\]    

Random ($\mathcal{N}$) infinitely divisible ($\mathcal{N}$ID) laws was introduced and developed by Gnedenko and Korolev (1996, section 4.6, p.137). This is based on $N_\theta$-sums, where $N_\theta$ is a positive integer-valued random variable (\textit{r.v}) having finite mean with probability generating function (\textit{p.g.f}) $P_\theta (s)= \varphi (\frac{1}{\theta}\varphi^{-1}(s)), \theta \in \Theta \subset (0,1)$. Here $\varphi$ is a Laplace transform (LT) which is also the standard solution to the Poincare equation, $\varphi(s) = P(\varphi(\theta s), s\geq 0, \theta \in \Theta$, $P$ being a p.g.f. The \textit{r.v} $N_\theta$ has the following property.

\pa {Lemma 1.2} $\theta N_\theta \stackrel{d}{\to}U$ as $\theta \downarrow 0$, and the LT of $U$ is $\varphi$, see  Gnedenko and Korolev (1996, p.138).

\vspace{4mm}
Now, we can study the limit distributions of random maximums using the transfer theorem for maximums by Gnedenko (1982). Satheesh \textit{et al.}(2008) briefly discussed the max-analogue of $\mathcal{N}$ID laws to obtain stationary solutions to a generalised max-AR(1) scheme. However, there was an  inadvertent omission, as the discussion did not stress that the LT $\varphi$ should also be the standard solution to the Poincare equation.

\vspace{4mm}
Proceeding from Satheesh \textit{et al.} (2008), we discuss random MID ($\mathcal{N}$MID) laws that is the max-analogue of $\mathcal{N}$ID laws, in section 2. In section 3 we discuss random max-stable laws and their domain of max-attraction, generalise certain results on GMS laws in Rachev and Resnick (1991) to random max-stable laws and the extremal processes associated with it.

\section{Random MID laws}
\setcounter{equation}{0}
\def\theequation{\thesection.\arabic{equation}} 

We begin by defining $\mathcal{N}$MID laws analogous to the $\mathcal{N}$ID laws of Gnedenko and Korolev (1996) correcting the omission mentioned above. Notice that we are describing a random maximum.

\pa{Definition 2.1} Let $\varphi$ be the standard solution to the Poincare equation and $N_\theta$ a positive integer-valued \textit{r.v} having finite mean with \textit{p.g.f} $P_\theta (s)= \varphi (\frac{1}{\theta}\varphi^{-1}(s)), \theta \in \Theta \subset (0,1)$. A \textit{d.f} $F(x)$ in $\mathbf{R^d}$ is $\mathcal{N}$MID if for each $\theta \in \Theta$ there exists a \textit{d.f} $G_\theta(x)$ that is independent of $N_\theta$, such that $F(x) = P_\theta (G_\theta (x))$ for all $x \in \mathbf{R^d}$.

\pa{Theorem 2.1} A \textit{d.f} $F(x)$ which is the limit of a sequence $F_n(x)$ of  $\mathcal{N}$MID \textit{d.f}s is itself $\mathcal{N}$MID.\\
\textit{Proof}. By virtue of the continuity of \textit{p.g.f}s, for every $\theta \in \Theta$ we have 
\[F(x)=\lim_{n\to\infty}F_n(x)= P_\theta (\lim_{n\to\infty}G_{\theta, n}(x))=P_\theta (G_{\theta}(x)).\]

Now we have an analogue of theorem 1.1, a de Finetti type theorem, for $\mathcal{N}$MID laws.
 
\pa{Theorem 2.2} Let $\varphi$ be the standard solution to the Poincare equation. A \textit{d.f} $F(x)$ in $\mathbf{R^d}$ is $\mathcal{N}$MID \textit{iff}  for some \textit{d.f}s $\{G_n\}$ and constants $\{a_n>0\}$

\[F(x)=\lim_{n\to\infty}\varphi\{a_n(1-G_n(x))\}.\]
\textit{Proof}. See the proof of theorem 3.5 in Satheesh \textit{et al.}(2008).

\vspace{4mm}
Notice that for a LT $\varphi(s), s>0$, $\varphi(\lambda (1-s)), 0<s\leq 1, \lambda>0$ is a p.g.f. Hence the above representation is essentially the weak limit of random-maximums under the transfer theorem for maximums. The next result facilitates the construction and/ or identification of $\mathcal{N}$MID \textit{d.f}s. 

\pa{Theorem 2.3} A \textit{d.f} $F(x)$ is $\mathcal{N}$MID \textit{iff} $F(x)=\varphi \{- \log H(x)\}$, where $\varphi$ is the standard solution to the Poincare equation and $H(x)$ a MID \textit{d.f}.\\
\textit{Proof}. We have seen that an $\mathcal{N}$MID \textit{d.f} $F(x)$ admits the representation for some \textit{d.f}s $G_\theta$  \[F(x)=\lim_{\theta \downarrow 0}\varphi\{\frac{1}{\theta}(1-G_\theta (x))\}.\]
Since $\varphi$ is continuous we can proceed as
 \[F(x)=\lim_{\theta \downarrow 0}\varphi \lo\{-\log \lo (\exp \{\frac{1}{\theta}(G_\theta (x) - 1)\}\ro) \ro\}= \varphi (-\log H(x))\]
 where $H(x) = \lim_{\theta \downarrow 0}  \exp \{\frac{1}{\theta}(G_\theta (x) - 1)\}$ is MID. Note the fact that every Poisson maximum is MID and every MID \textit{d.f} is the weak limit of Poisson maximums (Balkema and Resnick, 1977). Conversely, consider\\
\[ \varphi (-\log H(x))=  \int _0 ^\infty \exp \{t \log H(x)\} d\Lambda(t), t>0,\]
where $H(x)$ is MID and $\varphi$ is the LT of the \textit{d.f} $\Lambda.$ Now $\varphi (-\log H(x))$ is $\mathcal{N}$MID since the above is the integral representation of a \textit{d.f} that is the weak limit under the transfer theorem for maximums. This completes the proof.
 
\pa{Corollary 2.1} A \textit{d.f} is $\mathcal{N}$MID \textit {iff} it is the  limit distribution, as $\theta \downarrow 0$, of a random maximum of \textit{i.i.d r.v}s.

\vspace{4mm}    
We now proceed to prove the max-analogue of theorem 4.6.5 of Gnedenko and Korolev (1996, p.149). Let, for every $\theta \in \Theta,  \{X_{\theta,i}\}$ with \textit{d.f} $G_\theta$ be \textit{i.i.d} random vectors in in $\mathbf{R^d}$ and $N_\theta$ a positive integer-valued  \textit{r.v} having finite mean with \textit{p.g.f} $P_\theta (s)= \varphi (\frac{1}{\theta}\varphi^{-1}(s))$, that is independent of $\{X_{\theta,i}\}$, for every $\theta \in \Theta $ and $i$. Let $[\frac{1}{\theta}]$ denote the integer part of $\frac{1}{\theta}$. 

\pa{Theorem 2.4} Let $F(x)= \varphi(- \log G(x))$ be $\mathcal{N}$MID. Then
\be
\lim_{\theta \downarrow 0} P_\theta(G_\theta(x)) = \varphi(- \log G(x))  
\ee 
\textit{iff} there exists a \textit{d.f} $G(x)$ that is MID and 
\be
\lim_{\theta \downarrow 0} G_\theta ^{[\frac{1}{\theta}]}(x) = G(x).  
\ee
\textit{Proof}. The sufficiency of the condition (2.2) follows from the transfer theorem for maximums by invoking the relation  $\theta [\frac{1}{\theta}] \to 1$ and $\theta N_{\theta} \stackrel{d}{\to}U$ as $\theta \downarrow 0$ . Conversely (2.1) implies 
\be
\lim_{\theta \downarrow 0}\varphi (\frac{1}{\theta}\varphi^{-1}(G_{\theta}(x)))=\varphi(-\log G(x)).
\ee
Since $\varphi$ is a LT we have
\[\lim_{\theta \downarrow 0}\frac{1}{\theta}\varphi^{-1}(G_{\theta}(x))= -\log G(x)).\]
Again, since $\varphi (0)=1$, this implies that 
\be
\lim_{\theta \downarrow 0}G_{\theta}(x)=1.
\ee
Since $\varphi(\frac{1-G_{\theta}(x)}{\theta})$ is a \textit{d.f} that is $\mathcal{N}$MID for every $\theta \in \Theta$, by theorem 2.1, $\lim_{\theta \downarrow 0}\varphi(\frac{1-G_{\theta}(x)}{\theta})$ is also $\mathcal{N}$MID. Hence there exists a \textit{d.f} $H(x)$ that is MID such that 
\be
\lim_{\theta \downarrow 0} \lo(\frac{1-G_{\theta}(x)}{\theta}\ro) = -\log H(x).
\ee
On the other hand for $|\kappa| \leq 1$ we have
\be 
\log G_\theta ^{[\frac{1}{\theta}]} = \lo[\frac{1}{\theta}\ro] \log (1-(1-G_{\theta})) = \lo[\frac{1}{\theta}\ro] (G_{\theta}-1) + \kappa \lo[\frac{1}{\theta}\ro]|G_{\theta} - 1|^2.
\ee
Hence by (2.4) and (2.5) we get from (2.6)
\be
\lim_{\theta \downarrow 0} G_\theta ^{[\frac{1}{\theta}]}(x) = H(x).
\ee
Now applying the transfer theorem for maximums it follows that
\[\lim_{\theta \downarrow 0}\varphi(\frac{1}{\theta}\varphi^{-1}(G_{\theta}(x)))= \varphi(-\log H(x)).\]
Hence by (2.1) $H(x) \equiv G(x)$. That is, by (2.7), (2.2) is true with $G(x)$ being MID, completing the proof.

\section{Random max-stable laws}
\setcounter{equation}{0}
\def\theequation{\thesection.\arabic{equation}} 

Notice that the core of theorem 2.4 is that it identifies the limit of partial $N_\theta$-maximums of certain component \textit{r.v}s as a function of the limit of partial maximums of the same component \textit{r.v}s and vice-versa. This description thus enables us to define random max-stable ($\mathcal{N}$max-stable) laws analogous to the $\mathcal{N}$stable laws of Gnedenko and Korolev (1996) and their domains of $\mathcal{N}$max-attraction. This is facilitated by prescribing $\lo[\frac{1}{\theta}\ro]=n$ in theorem 2.4. Notice also that here the discussion can be for \textit{d.f}s in $\mathbf{R}$.

\pa{Definition 3.1} A \textit{d.f} $F(x)$ is $\mathcal{N}$max-stable \textit{iff} $F(x)=\varphi \{- \log H(x)\}$, where $\varphi$ is the standard solution to the Poincare equation and $H(x)$ a max-stable \textit{d.f}.

\pa{Theorem 3.1} An $\mathcal{N}$max-stable \textit{d.f} can be represented as $F(x) = P_\theta (F_\theta (x))$, for every $\theta \in \Theta$, where $F$ and $F_\theta$ are \textit{d.f}s of the same type. (This is in tune with definition 2.1).\\
\textit{Proof}. Since $F(x)$ is $\mathcal{N}$max-stable we have the following representation for every $\theta \in \Theta$
\[ F=\varphi \{- \log H\} = \varphi \lo( \frac{1}{\theta} \varphi ^{-1} (\varphi (-\theta \log H))\ro) = P_\theta (\varphi (-\theta \log H)) = P_\theta (\varphi (-\log H^\theta))=P_\theta(F_\theta).\]
Notice that $H$ and $H^\theta$ are \textit{d.f}s of the same type, Barakat \textit{et al.} (2009). Since $H$ is max-stable, $H^\theta$ also is max-stable. Thus the above representation describes an $\mathcal{N}$max-stable \textit{d.f} as an $N_\theta$-sum of \textit{d.f}s of the same type for every $\theta \in \Theta$, proving the result.

\vspace{4mm}
We now generalise Proposition 3.2 on GMS laws in Rachev and Resnick (1991) to $\mathcal{N}$max-stable laws.
\pa {Theorem 3.2} For a \textit{d.f} $F$ on $\mathbf{R}^d$ the following are equivalent.\\
(i) $F$ is $\mathcal{N}$max-stable\\
(ii) $\exp\{-\varphi^{-1}(F)\}$ is max-stable\\
(iii) There exists an $\ell \in [-\infty,\infty)^d$ and an exponent measure $\mu$ concentrated on $[\ell, \infty)$ such that for $x\geq \ell, F(x) = \varphi(\mu[\ell, x]^c)$.\\
(iv) There exists a multivariate extremal process $\{Y(t),t>0\}$ governed by a max-stable law and an independent \textit{r.v} $Z$ with \textit{d.f} $F$ and LT $\varphi$ such that $F(x) = P\{Y(Z) \leq x\}$.

\vspace{4mm}
\noindent \textit{Proof}. (i) $\Rightarrow$ (ii). $F$ is $\mathcal{N}$max-stable implies $F=\varphi\{-\log H\}$, where $H$ is max-stable. This implies $\exp\{-\varphi^{-1}(F)\} = H$ is max-stable.\\
(ii) $\Rightarrow$ (iii). From the representation of a max-stable \textit{d.f} by an exponent measure $\mu$ we have from (ii) $H(x) = \exp\{-\varphi^{-1} (F)\} = \exp\{-\mu([\ell,x]^c)]$.
This implies $\varphi^{-1} (F) = \mu([\ell,x]^c)$ or $F(x) = \varphi(\mu([\ell,x]^c).$\\
(iii) $\Rightarrow$ (iv). By (iii) we have an exponent measure $\mu$ corresponding to the max-stable law identified in (ii). Let $\{Y(t),t>0\}$ be the  extremal process governed by this max-stable law. That is 
\[ P\{Y(t)\leq x\} = \exp\{-t \mu([\ell,x]^c)\}.\]
Hence \[ P\{Y(Z)\leq x\}= \int _0 ^\infty \exp \{-t \mu([\ell,x]^c)\} dF(t) = \varphi(\mu([\ell,x]^c) = F(x).\]
(iv) $\Rightarrow$ (i) is now obvious. Thus the proof is complete.  

\vspace{4mm}
A notion that is closely associated with max-stable laws is their domain of max-attraction, Resnick (1987, p.12, 38, 263). The notion of geometric max-attraction for GMS laws were discussed by Rachev and Resnick (1991) and Mohan (1998). We now briefly discuss this for $\mathcal{N}$max-stable laws.

\pa{Definition 3.2} A \textit{d.f} $G(x)$ belongs to the domain of  $\mathcal{N}$max-attraction (D$\mathcal{N}$MA) of the \textit{d.f} $F(x)$ (with non-degenerate marginals) if there exists constants $a_{i,n}=a_i(\theta_n)>0$ and $b_{i,n}=b_i(\theta_n)$ such that 
\[\lim_{n \to \infty}P_n(G^n)=F,\]
meaning that  $\lim_{n \to \infty}P_n(G_i^n)=F_i$, for each $1\leq i \leq d$ where $G_i^n(x)=G_i^n(a_{i,n}x+b_{i,n})$ and $\theta_n=\frac{1}{n}.$
\vspace{4mm}   

Recalling that $\varphi$ is continuous and that max-attraction of $G$ to $H$ is equivalently specified by
\[n(1-G_i(a_{i,n}x+b_{i,n})) \rightarrow - \log H_i(x), 1\leq i \leq d,\] 
we have the following result as an immediate consequence of theorem 2.2. 

\pa{Theorem 3.3} Let $\varphi$ be the standard solution to the Poincare equation. A \textit{d.f} $F(x)=\varphi \{- \log H(x)\}$ is $\mathcal{N}$max-stable \textit{iff} for some \textit{d.f} $\{G\}$ and constants $a_{i,n}=a_i(\theta_n)>0$ and $b_{i,n}=b_i(\theta_n),$
\[\varphi(n(1-G_i(a_{i,n}x+b_{i,n}))) \rightarrow \varphi(- \log H_i(x))=F_i(x), 1\leq i \leq d.\] 

Again, from theorem 2.4, choosing $\theta$ such that $\lo[\frac{1}{\theta}\ro]=n$ and $G_\theta (x)= \lo(G_i(a_{i,n}x+b_{i,n}), 1\leq i \leq d \ro)$ where $a_{i,n}=a_i(\theta)>0$ and $b_{i,n}=b_i(\theta),$ from the classical results on max-stable laws and their domains of attraction, we have

\pa {Theorem 3.4}  A \textit{d.f} $G(x)$ belongs to the D$\mathcal{N}$MA of the \textit{d.f} $F(x)=\varphi \{- \log H(x)\}$ \textit{iff} it belongs to the DMA of $H(x)$.

\section{Conflict of Interests}
The authors declare that there is no conflict of interests regarding the publication of this article.


\begin{thebibliography}{99}
\bibitem{br}
{Balkema A A and Resnick}, 1977. Max-infinite divisibility. {J. App. Probab.}, {\bf14}, 309--319. 

\bibitem{bga}
{Barakat H M, Ghitany M E and Al-Hussaini E K,} 2009. Asymptotic distributions of order statistics and record values under the Marshall-Olkin parametrization operation. {Commu. Statist. - Theor. Meth.}, {\bf38}, 2267--2273.

\bibitem{gn}
{Gnedenko B V,} 1982. On limit theorems for a random number of random variables. In Fourth USSR-Japan Symp. Proc.. In Lecture Notes in mathematics, \textbf{Vol. 1021}, Springer, Berlin, 167--176.

\bibitem{gk}
{Gnedenko B V and Korolev V Yu.,} 1996. {Random Summation; Limit Theorems and Applications}, CRC Press, Boca Raton.

\bibitem{mo}
{Mohan N R,} 1998. On geometrically max infinitely divisible laws. {J. Ind. Statist. Assoc.}, {\bf36}, 01--12.

\bibitem{rr} 
{Rachev S T and Resnick S I,} 1991. Max-geometric infinite divisibility and stability. {Commu. Statist. - Stoch. Models}, {\bf7}, 191--218.
 
\bibitem{re}
{Resnick S I,} 1987. {Extreme Values, Regular Variation and Point Processes}. Springer-Verlag, New York.

\bibitem{ss}
{Satheesh, S, Sandhya E and Rajasekharan K E,} 2008. A generalisation and extension of an autoregressive model. {Statist. Probab. Lett.}, {\bf78}, 1369--1374. 
  
\end{thebibliography}
\end{document}